%% file: volume.tex
\documentclass[11pt]{amsart}

\usepackage{amsfonts, amstext, amsmath, amsthm, amscd, amssymb}
\usepackage{graphicx, color, epsfig, psfrag}

\newcommand{\tild}[1]{{\widetilde{#1}}}

\newcommand{\HH}{{\mathbb{H}}}

\newcommand{\tw}{{\rm tw}}

\theoremstyle{plain}
\newtheorem{theorem}{Theorem}[section]
\newtheorem{corollary}[theorem]{Corollary}
\newtheorem{lemma}[theorem]{Lemma}
\newtheorem{prop}[theorem]{Proposition}

\newtheorem*{no-num-theorem}{Theorem}

\theoremstyle{definition}
\newtheorem{define}[theorem]{Definition}

\begin{document}

\title{Volumes of highly twisted knots and links}
\author{Jessica S. Purcell}
\address{Jessica S. Purcell, Department of Mathematics, 1 University
  Station C1200, University of Texas at Austin, Austin, TX 78712}
\email{jpurcell@math.utexas.edu}

\begin{abstract}
  \input{abstract}

\end{abstract}


\maketitle


\section{Introduction \label{sec:intro}}
\input{introduction}

\section{Initial Geometric Estimates
  \label{sec:init-est}}
\input{init-est}

\section{Deformation Through Cone Manifolds \label{sec:deform-cone}}
\input{deform-cone}

\section{The Volume Estimate \label{sec:est-volume}}
\input{est-volume}

\section{Appendix \label{sec:appendix}}
\input{appendix}

\bibliographystyle{hamsplain}

\bibliography{references}


\end{document}

%% file: abstract.tex


We show that for a large class of knots and links with complements in
$S^3$ admitting a hyperbolic structure, we can determine bounds on the
volume of the link complement from combinatorial information given by
a link diagram.  Specifically, there is a universal constant C such
that if a knot or link admits a prime, twist reduced diagram with at
least 2 twist regions and at least C crossings per twist region, then
the link complement is hyperbolic with volume bounded below by 3.3515
times the number of twist regions in the diagram. C is at most 113.

%% file: introduction.tex

Given a diagram of a knot or link, our goal is to determine geometric
information about the complement of that link in $S^3$.  In
particular, if the complement admits a hyperbolic structure, then by
Mostow--Prasad rigidity that structure is unique.  We ought to be able
to make explicit statements about the geometry of this link
complement, including statements about its volume.  However, such
results based purely on a diagram seem to be rare.

Given a particular diagram, there are examples of volume computations.
For example, Cao and Meyerhoff proved the smallest volume knot or link
complement was the figure eight knot complement \cite{cao-meyerhoff}.
Given a simple diagram, computer software such as SnapPea, by Weeks,
can often compute a hyperbolic structure, including the volume
\cite{weeks:snappea}.

For particular classes of knots and links, other results on volume
have been determined.  Lackenby proved that in the special case in
which a knot or link is alternating, then the volume of the complement
is bounded above and below by the twist number of a diagram
\cite{lackenby:alt-volume}.  In fact, the upper bound is valid for all
knots, not just alternating.  This upper bound was further improved by
Agol and D. Thurston in an appendix to Lackenby's paper.
Additionally, they found a sequence of links with volume approaching
the upper bound.

Recently, the lower bound has been improved by results of Agol, Storm
and Thurston \cite{ast:volumes}.  The proof of this result still
requires that the links in question be alternating, however.

Our result is an extension of these results.  We are able to prove
that a similar lower bound holds for a large class of knots and links,
without the requirement that these knots be alternating.  We do,
however, need to introduce a requirement that the links be highly
twisted in each twist region.  Our methods use explicit deformation of
a hyperbolic structure through cone manifolds, using methods of
Hodgson and Kerckhoff \cite{hk:cone-rigid}, \cite{hk:univ}.  The high
amount of twisting ensures that the explicit estimates in their papers
and this one will apply throughout the deformation.

In order to state our results, we review some definitions concerning
the diagram of a knot or link, following Lackenby
\cite{lackenby:alt-volume}.

Consider the diagram of a knot or link as a 4--valent graph in the
plane, where associated to each vertex is over--under crossing
information.  A \emph{bigon} region is a region of this graph bounded
by only two edges.

\begin{figure}
  \begin{center}
    \includegraphics{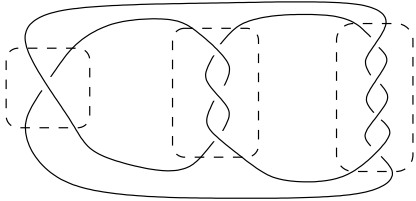}
  \end{center}
  \caption{Twist regions of the diagram are encircled with dashed
    lines.  For this diagram, $tw(D_K)=3$.}
  \label{fig:intro-twist-region}
\end{figure}

\begin{define}
  A \emph{twist region} of a diagram of a link $K$ consists of maximal
  collections of bigon regions arranged end to end.  A single crossing
  adjacent to no bigons is also a twist region.  Let $D_K$ denote the
  diagram of $K$.  We denote the number of twist regions in a diagram
  by $\tw(D_K)$.  See Figure \ref{fig:intro-twist-region}.
\end{define}

Our statements concern the number of twist regions of a diagram of a
knot or link.  In order to rule out extraneous twist regions, we need
the diagram to be reduced in the sense of the following two
definitions.

First, we require the diagram to be \emph{prime}.  That is, any simple
closed curve which meets two edges of the diagram transversely must
bound a region of the diagram with no crossings.

Second, we require the diagram to be \emph{twist reduced}.  That is,
if any simple closed curve meets the diagram transversely in four
edges, with two points of intersection adjacent to one crossing and
the other two adjacent to another crossing, then that simple closed
curve must bound a region of the diagram consisting of a (possibly
empty) collection of bigons arranged end to end between the crossings.

If a diagram of a hyperbolic knot or link is not prime, them some
crossings are extraneous and can be removed to obtain a prime diagram.
If it is not twist reduced, then a series of flypes will decrease the
number of twist regions of the diagram, until we are left with a twist
reduced diagram.

With these definitions, we are ready to state our results.

\begin{theorem}\label{thm:volumes}
  There is a universal constant $C$ such that if $K$ is a knot or link
  admitting a prime, twist reduced diagram $D_K$ with at least $2$
  twist regions and at least $C$ crossings in each twist region, then
  $S^3-K$ is hyperbolic with volume bounded below:
  $$\mbox{Volume}(S^3-K) \geq \tw(D_K)(3.3515).$$
  Further, the universal constant $C$ is at most $113$.
\end{theorem}

The upper bound on volume given by Lackenby, and Agol and D. Thurston
is also linear in the number of twist regions of the diagram
\cite{lackenby:alt-volume}.  Specifically, Agol and Thurston showed:
$$\mbox{Volume}(S^3-K) \leq 10\,v_3(\tw(D_K) -1)$$
where $v_3$ $(\approx 1.01494)$ is the volume of a regular hyperbolic
ideal tetrahedron.  

Thus the results in this paper extend the class of knots and links for
which volume is bounded above and below by linear functions of
$\tw(D_K)$.

Our proof is geometric in nature.  We begin with a link $L$ whose
complement is geometrically explicit.  We show that $S^3-K$ can be
obtained from $S^3-L$ by {D}ehn filling.  We determine a lower bound
on the volume of the link complement $S^3-L$.  We then perform
hyperbolic {D}ehn filling by a cone deformation.  Our final result is
obtained by bounding the change in geometry under this cone
deformation.

In Section \ref{sec:init-est}, we describe the link $L$ and explain
how to obtain $S^3-K$ from $S^3-L$ by {D}ehn filling.  In Section
\ref{sec:deform-cone}, we will review results on cone deformations and
describe their application to our particular situation.  Finally, in
Section \ref{sec:est-volume} we put these results together to conclude
the proof of Theorem \ref{thm:volumes}.

%% file: init-est.tex

\subsection{The augmented link \boldmath{$L$}}
  
Start with a prime, twist reduced diagram $D_K$ of a link $K$.  In the
rest of this paper, we will assume the diagram $D_K$ is fixed once and
for all.  To simplify notation, we will refer to this fixed diagram
and the knot by the same symbol $K$.  Thus we will write $\tw(K)$ to
mean $\tw(D_K)$.  The reader should note that the twist number
$\tw(D_K)$ \emph{is} dependent upon the particular choice of prime,
twist reduced diagram.  However, since our diagram is now fixed, this
simplification should cause no confusion.  

\begin{figure}
  \begin{center}\includegraphics{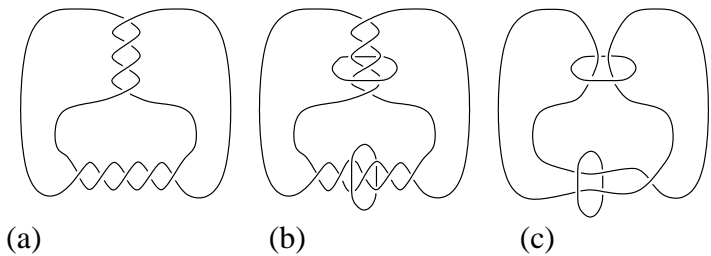}\end{center}
  \caption{(a) The original link diagram $K$.  (b) Crossing circles
    added.  (c) The diagram of $L$.}
  \label{fig:knot-decomp}
\end{figure}

Given our diagram of $K$, we obtain a new link by adding additional
link components to the diagram.  At each twist region, encircle the
twist region by a simple curve, called a \emph{crossing circle}.  See
Figure \ref{fig:knot-decomp} (a) and (b).  Links with added crossing
circles have been studied by many people, including Adams
\cite{adams:aug}.  These links were used by Lackenby, and also by Agol
and D. Thurston to improve Lackenby's volume results for alternating
links \cite{lackenby:alt-volume}.  Provided the original diagram of
$K$ was prime and twist reduced with at least two twist regions, then
the link with crossing circles added is known to be hyperbolic.  This
can be shown either using methods of Adams \cite{adams:aug}, or
directly using {A}ndreev's theorem \cite{purcell:thesis}.

Let $J$ be the link with crossing circles added to each twist region
of $K$.  Modify the diagram of $J$ by removing pairs of crossings at
each twist region, and let $L$ be this new link.  Now, $S^3-J$ is
homeomorphic, and thus isometric by Mostow--Prasad rigidity, to the
manifold $S^3-L$.  The diagram of $L$ consists of strands in the
projection plane encircled in pairs by crossing circles.  These
strands in the projection plane only cross, if at all, in pairs at
crossing circles, and here they may only cross once.  We will call
such a link an \emph{augmented link}.  See Figure
\ref{fig:knot-decomp} (c).

The augmented link $L$ has a complement with nice geometric
properties, allowing us to give estimates on the volume of $S^3-L$.
Before giving these estimates, we relate $S^3-L$ to the complement of
the original link, $S^3-K$.

\subsection{{D}ehn filling on \boldmath{$S^3-L$}}

Let $C_i$ be the $i$-th crossing circle of the link $L$.  Take a
horoball neighborhood $N_i$ of the cusp corresponding to $C_i$ in
$S^3-L$. Then $N_i$ has torus boundary, with meridian $\mu$ and
longitude $\lambda$.  If we glue a solid torus onto $S^3-L-N_i$ in
such a way that the boundary of the solid torus is glued to the
boundary $\partial N_i$, with the curve $\mu + n_i\lambda$ on
$\partial N_i$ bounding a meridional disk of the solid torus, then we
obtain a new link complement in $S^3$.  The diagram of this link is
the same as that of $L$, only now $n_i$ full twists (i.e. $2n_i$
crossings) have been inserted into the $i$-th twist region, and the
crossing circle $C_i$ has been removed (see for example
\cite{rolfsen-book}, Chapter 9).  Recall that the insertion of a solid
torus into the manifold $S^3-L$ in this manner is called the
\emph{{D}ehn filling of $S^3-L$ along the slope $1/n_i$ on $C_i$}.

Now, since $L$ was originally obtained from $K$ by adding crossing
circles and removing crossings at twist regions in a diagram, by
choosing values for the $n_i$ appropriately at each crossing circle,
this type of {D}ehn filling performed at each crossing circle will
give us back the original link complement $S^3-K$.  Precisely, we
perform {D}ehn filling of $S^3-L$ along the slopes $(1/n_1, 1/n_2,
\dots, 1/n_{\tw(K)})$ on crossing circle components to obtain the
manifold $S^3-K$.

\subsection{Decomposition of \boldmath{$S^3-L$}}

In their appendix to Lackenby's paper \cite{lackenby:alt-volume}, Agol
and D. Thurston describe how to decompose the link complement $S^3-L$
into totally geodesic ideal polyhedra when $L$ happens to have no
crossings of strands in the projection planes (i.e., in our case, when
$K$ had an even number of crossings at each twist region).  Their
methods immediately extend to the case when single crossings of the
strands in the projection plane are allowed at each twist region.  In
the following paragraphs, we review the main results.  See also
\cite{futer-purcell} for more details and pictures. 

First, consider a link with no single crossings in the projection
plane.  Call this link $\bar{L}$.  Notice that a reflection through
the projection plane preserves $\bar{L}$.  Thus the projection plane
is a totally geodesic surface in $S^3-\bar{L}$ (see for example
\cite{leininger}).  Now slice $S^3-\bar{L}$ along the projection
plane.  This breaks the manifold into two pieces with totally geodesic
boundary.  Each 2--punctured disc bounded by a crossing circle has
been sliced in half.  Next slice each of these halves and open them up
into two triangles.  Since triangles can be taken to be totally
geodesic, the result is two identical polyhedra, $P_1$ and $P_2$, with
totally geodesic faces.

\begin{figure}
  \begin{center}
    \includegraphics[width=1.4in]{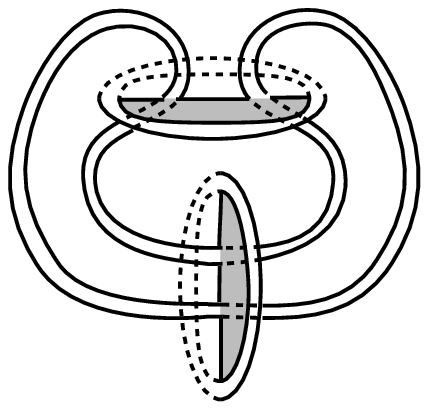}
    \hspace{.1in}
    \includegraphics[width=1.4in]{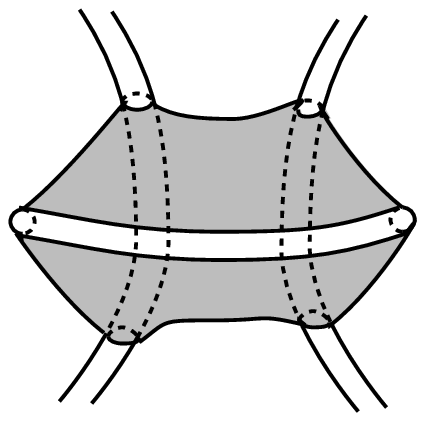}
    \hspace{.1in}
    \includegraphics[width=1.5in]{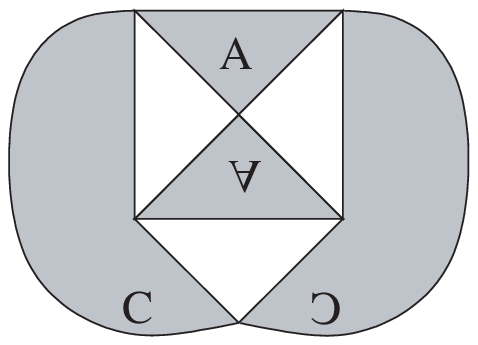}
  \end{center}
  \caption[Decomposing $S^3 - L$ into ideal
  polyhedra.]{Decomposing $S^3 - L$ into ideal polyhedra:
    First slice along the projection plane, then split remaining
    halves of two--punctured disks.  Obtain polyhedron on right.}
\label{fig:aug-decomp}
\end{figure}

The ideal polyhedra $P_1$ and $P_2$ have two kinds of faces.  ``White
faces'' come from regions of the projection plane.  The triangular
``shaded faces'' consist of halves of the 2--punctured discs bounded
by crossing circles.  These two types of faces intersect each other at
right angles.  The intersections are the edges of the polyhedra.  The
edges meet in 4--valent ideal vertices.  At these vertices, two white
and two shaded faces lie across from each other in pairs.  See Figure
\ref{fig:aug-decomp}.

We obtain the manifold $S^3-\bar{L}$ from these polyhedra by reversing
the slicing procedure above.  First, on each $P_i$, form the
2--punctured disc halves by gluing back together the two shaded
triangles which were opened along their common vertex.  Then glue the
white faces of $P_1$ to the identical white faces of $P_2$.

By changing the gluing procedure slightly, we can obtain the manifold
$S^3-L$ from the polyhedra $P_1$ and $P_2$ as well.  Recall the link
$L$ may have single crossings at some crossing circles.  For each
crossing circle $C_i$ of $L$, consider the four shaded triangles (two
on $P_1$, and two on $P_2$) making up the disc bounded by the
corresponding crossing circle of $\bar{L}$.  For example, the shaded
faces labelled with A in Figure \ref{fig:aug-decomp} are two of
the four triangles.  The other two lie on an identical polyhedron.

\begin{figure}
  \begin{center}
    \includegraphics{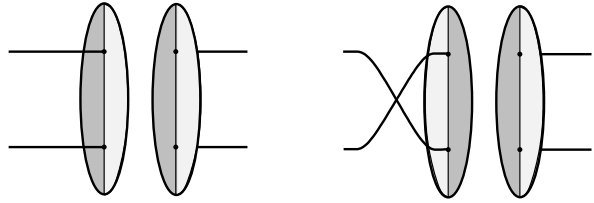}
  \end{center}
  \caption{Left:  Gluing 2--punctured discs with no crossings.  Right:
    Gluing 2--punctured discs with a single crossing.}
  \label{fig:half-twist}
\end{figure}

If there are no single crossings at the twist region of $C_i$, then
the gluing is the same as for $S^3-\bar{L}$: Glue the triangles of
$P_i$, $i=1, 2$, to each other across their common vertex.  For
example, in Figure \ref{fig:aug-decomp}, the faces labelled A
would be glued together.  If there is a single crossing at $C_i$, glue
each triangle of $P_1$ to the opposite one of $P_2$, matching vertices
from the crossing circle.  For example, if we replace the link on the
left in Figure \ref{fig:aug-decomp} with one with a single crossing in
the top twist region, then in the gluing of the two polyhedra the top
face labelled A in the figure on the right would be glued to the
opposite face labelled A of the opposite polyhedron.  This puts a
``half--twist'' into the manifold.  See Figure \ref{fig:half-twist}
for another schematic picture.  In either case, $S^3-L$ is then
obtained by gluing corresponding white faces of $P_1$ and $P_2$
together, as in the gluing of $S^3-\bar{L}$.

Now, we wish to find lower bounds on volume.  We will give a rough
lower bound first by considering the cusp volume of $S^3-L$.

\subsection{Volume estimate for \boldmath{$S^3-L$}}

To give an estimate on the cusp volume of $S^3-L$, we must analyze
neighborhoods of the ideal vertices of the polyhedra.  These
neighborhoods were analyzed quite carefully with Futer in
\cite{futer-purcell}.  Again we record the results here.

A horospherical torus about a cusp of the manifold $S^3-L$ will
intersect the polyhedra $P_1$ and $P_2$ in rectangles.  This is
because at a cusp, ideal vertices of $P_1$ and $P_2$ meet. These
vertices are all 4--valent, meaning their intersection with a
horosphere is a quadrilateral, and all faces of the $P_i$ meet at
right angles, forcing the quadrilateral to be a rectangle.  The sides
of this rectangle can be colored according to the color of the faces
of intersection.  Thus each rectangle will have two ``white sides''
across from each other, and two ``shaded sides'' across from each
other.  When we glue the polyhedra together to form $S^3-L$, we glue
these rectangles along white and shaded sides.  Thus each cusp of
$S^3-L$ is tiled by rectangles.  The number and size of these
rectangles will give us a volume estimate.

The following two lemmas are essentially Lemmas 2.3 and 2.6 of
\cite{futer-purcell}.

\begin{lemma}
  A horospherical torus about a cusp in $S^3-L$ corresponding to a
  crossing circle is tiled by two rectangles.  These rectangles are
  given by the intersection of the horospherical torus with $P_1$ and
  $P_2$.

\begin{itemize}
\item  A longitude of the crossing circle is isotopic to the curve given by
  stepping along two shaded sides of the rectangles.
  
\item  A meridian is isotopic to the curve given by a white side of a
  rectangle plus $\epsilon$ times a shaded side, where $\epsilon=0,\pm
  1$ depending on whether the crossing circle bounded no crossings, or
  one crossing in the positive or negative direction, respectively.
  \label{lemma:cross-circ}
\end{itemize}
\end{lemma}

\begin{lemma}
  A horospherical torus about a cusp in $S^3-L$ corresponding to a
  link component $L_i$ in the projection plane is tiled by $2\,C(L_i)$
  rectangles, where $C(L_i)$ is the number of crossing circles,
  counted with multiplicity, through which $L_i$ passes in the diagram
  of $L$.
  \label{lemma:knot-strand-1}
\end{lemma}

When $K$ is a knot with diagram $D_K$, there is just one link
component in the projection plane, and so it passes through each
crossing circle twice.  Since the number of crossing circles in the
diagram of $L$ is $\tw(K)$, in this case the number of rectangles in
the cusp of $S^3-L$ corresponding to the link component in the
projection plane is $4\,\tw(K)$.

When $K$ is a link, we still have two strands passing through each
crossing circle.  Thus the total number of rectangles in all the cusps
of $S^3-L$ corresponding to link components on the projection plane is
still $4\,\tw(K)$.

We record these results in the following lemma.

\begin{lemma}
  There are a total of $4\,\tw(K)$ rectangles tiling all
  horospherical tori about the cusps of $S^3-L$ corresponding to link
  components in the projection plane.
  \label{lemma:knot-strand}
\end{lemma}

We need the sizes of these rectangles.  Again these results were
previously computed.  The following is essentially Corollary 3.9 of
\cite{futer-purcell}.

\begin{lemma}
  There exists a choice of horospheres expanded about all cusps of
  $S^3-L$ so that given this expansion, the length of any shaded side
  of a rectangle is exactly one.  The length of any white side of a
  rectangle is at least one.
\label{lemma:rect_lengths}
\end{lemma}

These lemmas immediately give the following proposition.

\begin{prop}
  Let $K$ be a link with a prime, twist reduced diagram with
  $\tw(K)\geq 2$ twist regions.  Then the augmented link $L$ obtained
  by adding a crossing circle to the diagram of $K$ at each twist
  region has complement $S^3-L$ with cusp volume at least $3\,\tw(K)$.
  \label{prop:init-vol}
\end{prop}

\begin{proof}
  We can compute the cusp volume by computing the volume of a region
  in the universal cover $\HH^3$ of $S^3-L$ which projects to the cusp
  in a one--to--one manner.
  
  Consider the upper half space model of $\HH^3$.  For any cusp, we
  can conjugate so that the point at infinity in this model of $\HH^3$
  projects to that cusp under the covering map.  Then the rectangles
  tiling that cusp lift to give rectangles tiling a horosphere about
  infinity.  For each rectangle on the cusp, choose a representative
  on the horosphere projecting to that rectangle.  Then the volume of
  the cusp is given by the sum of the volumes of regions of $\HH^3$
  lying over these rectangles.

  By Lemma \ref{lemma:rect_lengths}, the volume over any of these
  rectangles is at least as large as the volume lying over the square
  with side length $1$.  A calculation shows that the volume lying
  over a square of side length $1$ is $1/2$.  So the volume is at
  least $1/2$ times the total number of rectangles tiling all cusps of
  $S^3-L$.
  
  To determine the total number of rectangles tiling all cusps of
  $S^3-L$, we use the previous lemmas.  By Lemma
  \ref{lemma:knot-strand}, there are $4\,\tw(K)$ rectangles tiling
  all cusps coming from link components in the projection plane.  By
  Lemma \ref{lemma:cross-circ}, there are $2$ rectangles per crossing
  circle.  But there are $\tw(K)$ total crossing circles in the
  diagram of $L$.  Hence there are $2\,\tw(K)$ rectangles tiling
  cusps coming from crossing circles.  So the total number of
  rectangles tiling all cusps of $S^3-L$ is $6\,\tw(K)$.
  
  Thus the total cusp volume is at least $6\,\tw(K)(1/2) =
  3\,\tw(K)$.
\end{proof}

We obtain a better estimate on volume by combining this cusp estimate
with B{\"o}r{\"o}czky's lower bound for the density of a horoball
packing in hyperbolic space \cite{boroczky}.  By his result, the
volume of a maximal cusp neighborhood in a hyperbolic 3--manifold $M$
is at most $\sqrt{3}/(2v_0) \approx 0.853276$ times the volume
of the manifold $M$, where $v_0$ is the volume of a regular ideal
hyperbolic tetrahedron.  Thus the volume of $S^3-L$ is at least
$(3\,\tw(K))/(0.853276) \approx 3.51586\, \tw(K)$.

\begin{corollary}
  Let $K$ be a link with a prime, twist reduced diagram with
  $\tw(K)\geq 2$ twist regions.  Then the augmented link $L$ obtained
  by adding a crossing circle to the diagram of $K$ at each twist
  region has complement $S^3-L$ with volume at least $(3.51586)\,\tw(K)$.
  \label{cor:init-vol}
\end{corollary}

%% file: deform-cone.tex

\subsection{Hyperbolic {D}ehn filling}

In the previous section, we obtained a hyperbolic link whose geometry
could be determined explicitly.  This hyperbolic link was related to
our original link by a {D}ehn filling, i.e., by gluing solid tori into
the cusps corresponding to crossing circles.

Thurston showed that most {D}ehn fillings on components of a
hyperbolic link complement can be obtained by a hyperbolic {D}ehn
filling, that is, by deforming the complete hyperbolic geometric
structure on the cusped manifold through incomplete hyperbolic
structures, yielding a new manifold with a final complete hyperbolic
structure \cite{thurston}.  Hodgson and Kerckhoff were able to make
this explicit in \cite{hk:univ}.  They showed that if the normalized
lengths of the slopes along which the hyperbolic {D}ehn filling is
performed are each longer than a universal constant, then hyperbolic
{D}ehn filling is possible.  In particular, in this case of long
slopes, a special type of deformation through incomplete hyperbolic
structures exists, called a hyperbolic cone deformation.  In a
hyperbolic cone deformation, each of the intermediate incomplete
structures of the deformation is a \emph{hyperbolic cone manifold}.

For a complete description of hyperbolic cone manifolds and the
metrics involved, see \cite{hk:cone-rigid}.  For our purposes, a
hyperbolic cone manifold $M$ admits a smooth hyperbolic metric
everywhere except along a link $\Sigma$.  At each link component
$\Sigma_i$ of $\Sigma$, a tubular neighborhood has meridional cross
section which is a 2--dimensional hyperbolic cone, with cone angle
$\alpha_i$.  The cone manifold structure on $M$ is locally
parameterized by the collection of angles $(\alpha_1, \alpha_2, \dots,
\alpha_N)$.

In our case, recall that we wish to perform {D}ehn filling on $S^3-L$
to obtain $S^3-K$.  The fillings are performed along the slopes
$1/n_i$ on the crossing circles.  Thus our singular locus $\Sigma$
will consist of the cores of the crossing circles.  Provided our
slopes are long enough to meet the requirements found in
\cite{hk:univ}, we will obtain a cone deformation of the hyperbolic
structure on $S^3-L$ with initial cone angles all $0$.  Each cone
angle will increase strictly monotonically to cone angle $2\pi$ under
the deformation.

In \cite{hk:univ}, Hodgson and Kerckhoff found bounds on the change of
volume under certain cone deformations, in particular those with a
single component of the singular locus $\Sigma$.  Since we are filling
along multiple cusps, with multiple components of $\Sigma$ (one per
crossing circle), we will need versions of the results in
\cite{hk:univ} which allow filling along more than one cusp.  We state
those results here.

\subsection{Change of Volume}

The following theorem is Theorem 6.5 of \cite{hk:univ}, modified for
the case of multiple components of the singular locus $\Sigma$.  The
notation is from that paper.  In particular,
\begin{equation}
  H(z) = \frac{1+z^2}{3.3957\,z(1-z^2)},
  \label{eqn:H}
\end{equation}
and
\begin{equation}
  \tild{G}(z) = \frac{(1+z^2)^2}{6.7914\,z^3(3-z^2)}.
  \label{eqn:Gtild}
\end{equation}

\begin{theorem}
  Let $X$ be a cusped hyperbolic 3--manifold and $M$ a hyperbolic
  3--manifold which can be joined by a smooth family of hyperbolic
  cone manifolds with cone angles $0\leq\alpha_i\leq 2\pi$ along the
  $i$-th component $\Sigma_i$ of a link $\Sigma$.  Suppose that
  $\alpha_i\ell_i \leq 0.5098$ holds throughout the deformation,
  where $\ell_i$ denotes the length of $\Sigma_i$.  Then the
  difference in volume
  $$ \Delta V = \mbox{Volume}(X) - \mbox{Volume}(M)$$
  satisfies
  $$
  \Delta V \leq \sum_i \int_{\hat{z}_i}^1
  \frac{H'(w)}{8H(w)(H(w)-\tild{G}(w))}\,dw.
  $$
  Here $\hat{z}_i = \tanh(\hat{\rho}_i)$, $\hat{\rho}_i$ is the unique
  solution of $1.69785 \tanh(\hat{\rho}_i)/\cosh(2\hat{\rho}_i) =
  2\pi\hat{\ell}_i$ with 
  $\hat{\rho}_i \geq 0.531$, and $\hat{\ell}_i$ is the length of
  $\Sigma_i$ in $M$.
  \label{thm:hk-volumes}
\end{theorem}

Since this version of the theorem is somewhat different from that in
\cite{hk:univ}, we outline its proof in the appendix, Section
\ref{sec:appendix}.

%% file: est-volume.tex

In this section, we will find conditions under which Theorem
\ref{thm:hk-volumes} will apply.  Using these conditions, we can
complete the proof of Theorem \ref{thm:volumes}.

For Theorem \ref{thm:hk-volumes} to apply, we first need a smooth cone
deformation from the manifold $S^3-L$ to the manifold $S^3-K$.  We
also need each pair of cone angle $\alpha_i$ and length of singular
locus component $\ell_i$ to satisfy the inequality $\alpha_i\ell_i
\leq 0.5098$.

In fact, the existence of the deformation and the bound on
$\alpha_i\ell_i$ are both given by Hodgson and Kerckhoff in Theorem
5.12 of \cite{hk:univ}. (More accurately, these are given in the proof
of that theorem.)  In our situation, this theorem states that provided
the slopes $1/n_i$ along which we perform {D}ehn filling have
normalized length at least $\sqrt{2(56.4696)} \approx 10.628$, then a
cone deformation exists with initial point the complete hyperbolic
structure on $S^3-L$, and extends to give the complete hyperbolic
structure on $S^3-K$ with each cone angle equal to $2\pi$.

Normalized length is defined as follows.

\begin{define}
  Let $s$ be a slope on a cusp with cusp torus $T$.  Let $\alpha$ be a
  geodesic representative of $s$ on $T$.  The \emph{normalized length}
  of $s$ is defined to be the length of $\alpha$ divided by the square
  root of the area of $T$:
  $$\mbox{Normalized length}(s) =
  \frac{\mbox{Length}(\alpha)}{\sqrt{\mbox{Area}(T)}}.$$
\end{define}

At the $i$-th crossing circle, recall we perform {D}ehn filling along
the curve $\sigma_i = \mu + n_i \lambda$, where $\mu$ is a meridian of
a torus about this crossing circle, and $\lambda$ is a longitude.
Thus we need to estimate the normalized length of this curve, to
ensure it is at least $\sqrt{2(56.4696)}$.

In Section \ref{sec:init-est}, Lemma \ref{lemma:cross-circ}, we saw
that a torus about a crossing circle is tiled by two rectangles with
white and shaded sides, and meridian and longitude given by steps
along those sides.  The longitude $\lambda$ is given by two steps
along shaded sides.  If the crossing circle bounds a single crossing,
then $\mu$ is given by one step along a white side, plus or minus a
step along a shaded side, and $\sigma_i$ is given by a step along a
white side plus (or minus) $2n_i + 1$ steps along shaded sides.  If
the crossing circle bounds no single crossing, then $\mu$ is given by
one step along a white side, and $\sigma_i$ is given by a step along a
white side plus (or minus) $2n_i$ steps along shaded sides.  In either
case, note that if $c_i$ is the number of crossings in the $i$-th
twist region, then the curve the curve $\sigma_i$ is given by one step
along a white side plus (or minus) $c_i$ steps along shaded sides.

Let $w$ denote the length of a white side and let $s$ denote the
length of a shaded side.  Since the area of the torus is $2sw$, the
normalized length $\widehat{L}_i$ of $\sigma_i$ is given by:
$$\widehat{L}_i = \sqrt{\frac{w}{2s}+\frac{c_i^2s}{2w}}.$$
This will be
a minimum when $w/2s$ equals $c_i/2$, and that minimum value is
$\sqrt{c_i}$.  Hence we need $c_i$ large enough that
$$\sqrt{c_i} \geq \sqrt{2(56.4696)} = \sqrt{112.9392}$$
Thus $c_i \geq 113$ will be sufficient.

Then Theorem 5.12 of \cite{hk:univ} will apply to $X=S^3-L$ to give
the desired cone deformation and bounds on the $\alpha_i\ell_i$.
Hence Theorem \ref{thm:hk-volumes} applies to $X=S^3-L$ and $M=S^3-K$,
and we find:
$$\mbox{Volume}(X) - \mbox{Volume}(M) \leq \sum_{i=1}^{\tw(K)}
\int_{\hat{z}_i}^1 \frac{H'(w)dw}{8H(w)(H(w)-\tild{G}(w))}$$
Note the sum is over all components of the singular locus.  In our
case, these correspond to the crossing circles, and there are exactly
$\tw(K)$ of these. 

In the $i$-th cusp, the value of $z$ decreases from $1$ to the value
$\hat{z}_i$.  Recall from the statement of Theorem
\ref{thm:hk-volumes} that $\hat{z}_i = \tanh(\hat{\rho}_i)$, where
$\hat{\rho}_i$ is the unique solution of some equation such that
$\hat{\rho}_i \geq 0.531$.  In particular, $\hat{\rho}_i$ is
guaranteed to be at least $0.531$.  So $\hat{z}_i$ is at least $z_1 =
\tanh(0.531)$.  Thus the decrease in volume is at most:
\begin{equation*}
\sum_{i=1}^{\tw(K)} \int_{z_1}^1
\frac{H'(w)dw}{8H(w)(H(w)-\tild{G}(w))} \leq
\sum_{i=1}^{\tw(K)} 0.16436 = 
\tw(K)(0.16436)
\end{equation*}

Hence
$$\mbox{Volume}(M) \geq \mbox{Volume(X)} - \tw(K)(0.16436)$$

Also, we know by Corollary \ref{cor:init-vol} that $\mbox{Volume}(X)
\geq (3.51586)\,\tw(K)$.  Hence we obtain our final result:
$$\mbox{Volume}(M) \geq (3.51586)\,\tw(K) - \tw(K)(0.16436) =
\tw(K)(3.3515)$$

This concludes the proof of Theorem \ref{thm:volumes} which we
restate.

\begin{theorem}
  Let $K$ be a link in $S^3$ admitting a prime, twist reduced diagram
  with $\tw(K)\geq 2$ twist regions and at least $113$ crossings per
  twist region.  Then the volume of $S^3-K$ satisfies:
  $$\mbox{Volume}(S^3-K) \geq \tw(K)(3.3515)$$
\end{theorem}

%% file: appendix.tex

In this appendix, we sketch the proof of Theorem
\ref{thm:hk-volumes}.  

The proof is identical to that in \cite{hk:univ}, except for small
modifications that need to be made for the case in which the cone
manifold has multiple components of the singular locus.

In particular, the following changes to that proof need to be made.
First, in \cite{hk:univ}, the cone deformation was parameterized by
$t=\alpha^2$.  When there are multiple components, we can no longer
guarantee that the deformation is parameterized by $t=\alpha_j^2$ for
each $\alpha_j$.  We do know that some parameterization exists for
which the calculations of \cite{hk:univ} go through (see the comment
at the bottom of page 41 in \cite{hk:univ}, and Purcell
\cite{purcell:cusps}), but we don't know specifically what that
parameterization will be.  Thus we need to modify the calculations to
be independent of parameterization.

We also need to make modifications to calculations which used the area
estimate for a single cusp, Theorem 4.4 of \cite{hk:univ}.  As stated
in that theorem, the lower bound for the area is half as large when
the singular locus has multiple components.  Thus we replace Corollary
5.1 of \cite{hk:univ} with half that estimate:
\begin{equation}
\alpha_i\ell_i \geq \frac{1}{2}\left(
  3.3957\frac{\tanh(R)}{\cosh(2R)}\right) =
1.6978\frac{\tanh(R)}{\cosh(2R)}.
\label{eqn:h}
\end{equation}
We recompute calculations using this new estimate.

With these two modifications, the proof goes through nearly as written
in \cite{hk:univ}.  Following their notation, we begin by letting $u_i
= \alpha_i/\ell_i$.  Hodgson and Kerckhoff bounded the change in $u_i$
by finding bounds on $du_i/dt$, in Proposition 5.6.  Rather than take
derivatives with respect to time $t$, we take derivatives with respect
to cone angle $\alpha_i$ in the $i$-th cusp.  This is independent of
parameterization.  The half of Proposition 5.6 of \cite{hk:univ}
necessary for Theorem \ref{thm:hk-volumes} becomes:
\begin{equation}
\frac{1}{\alpha_i}\frac{du_i}{d\alpha_i} \leq 4\tild{G}(z_i)
\label{eqn:prop5.6}
\end{equation}
where
the function $\tild{G}(z)$ is defined as in (\ref{eqn:Gtild}).

Now we are ready to modify the proof of Theorem 6.5 of \cite{hk:univ}
directly.  From the {S}chl{\"a}fli formula (see
\cite{Cooper-HK:orbifolds}, \cite{hodgson:thesis}), the change in
volume $V$ of a hyperbolic cone manifold during a deformation
satisfies
$$dV = -\frac{1}{2}\sum_i \ell_i \, d\alpha_i = -\sum_i
\frac{\alpha_i\,d\alpha_i}{2u_i}.$$

We will rewrite $u_i$ in terms of the function $H$, where $H$ is given
by (\ref{eqn:H}).  $H$ was defined so that $H(\tanh(R))$ is the
reciprocal of the function
$$h(R) = 3.3957 \frac{\tanh(R)}{\cosh(2R)},$$
which appears in the
inequality (\ref{eqn:h}).  Provided $R \geq 0.531$, this function
$h(R)$ is strictly decreasing, and thus we may define its inverse for
$x$ in the interval $(0, h(0.531)] \approx (0, 1.019675]$.  Since
$\alpha_i\ell_i \leq 0.5098$, $2\alpha_i\ell_i \leq h(0.531)$.  We
define $\rho_i = h^{-1}(2\alpha_i\ell_i)$.  Letting $z_i =
\tanh(\rho_i)$, we have $u_i = 2\alpha_i^2H(z_i)$.  So the term
$-\alpha_i/(2u_i)$ equals $-1/(4\alpha_i H(z_i))$, and we obtain:
\begin{equation}
  dV = -\sum_i \frac{1}{4\alpha_i H(z_i)} d\alpha_i
  \label{eqn:dV}
\end{equation}

Now, by inequality (\ref{eqn:prop5.6}), we have:
$$\frac{1}{\alpha_i}\frac{du_i}{d\alpha_i} \leq 4
\widetilde{G}(z_i).$$
We compute:
$$\frac{1}{\alpha_i}\frac{du_i}{d\alpha_i} =
\frac{1}{\alpha_i}\frac{d}{d\alpha_i}(2\alpha_i^2H(z_i)) =
4H(z_i)+2\alpha_iH'(z_i)\frac{dz_i}{d\alpha_i}.$$

Hence
$$\alpha_iH'(z_i)\frac{dz_i}{d\alpha_i} \leq
2(\widetilde{G}(z_i)-H(z_i))$$
which implies
\begin{equation}
  -\frac{1}{\alpha_i} \geq
\frac{H'(z_i)}{2(H(z_i)-\widetilde{G}(z_i))} \frac{dz_i}{d\alpha_i}.
\label{eqn:alpha_ineq}
\end{equation}
To see that the direction of the inequalities is correct, notice that
$H(z_i)-\widetilde{G}(z_i)$ is positive.  

Now, multiply the inequality (\ref{eqn:alpha_ineq}) by $1/(4H(z_i))$.
This bounds each term in equation (\ref{eqn:dV}).
$$dV \geq \sum_i
\frac{H'(z_i)}{8H(z_i)(H(z_i)-\widetilde{G}(z_i))}dz_i.$$

We integrate over the deformation.  As $\alpha_i$ increases from $0$
to $2\pi$, $z_i$ decreases from $1$ to $\hat{z}_i$.  When $\alpha_i$
has reached $2\pi$, it remains $2\pi$ and there is no further
contribution.

Hence
$$ \Delta V \leq
\sum_i\int_{\hat{z}_i}^1
\frac{H'(z_i)}{8H(z_i)(H(z_i)-\widetilde{G}(z_i))}dz_i.$$